\documentclass{article}%{amsart}
\usepackage{amssymb,amsmath}%,showkeys}
\usepackage[mathscr]{euscript}
\input{epsf}

\newcounter{sec}

\newcounter{punct}[sec]

\def\punct{\refstepcounter{punct}{\arabic{sec}.\arabic{punct}.  }}

\newtheorem{theorem}{Theorem}[sec]

\newtheorem{lemma}[theorem]{Lemma}

\newtheorem{corollary}[theorem]{Corollary}

\def\COUNTERS{\addtocounter{sec}{1}
              \setcounter{punct}{0}
          \setcounter{equation}{0}
          \setcounter{theorem}{0}
          }
          
          \def\sm{\smallskip}

\begin{document}

\newcommand{\supp}{\mathop {\mathrm {supp}}\nolimits}
\newcommand{\rk}{\mathop {\mathrm {rk}}\nolimits}
\newcommand{\Aut}{\mathop {\mathrm {Aut}}\nolimits}
\newcommand{\Out}{\mathop {\mathrm {Out}}\nolimits}
\renewcommand{\Re}{\mathop {\mathrm {Re}}\nolimits}

\def\ov{\overline}

\def\wt{\widetilde}

\renewcommand{\rk}{\mathop {\mathrm {rk}}\nolimits}
\renewcommand{\Aut}{\mathop {\mathrm {Aut}}\nolimits}
\renewcommand{\Re}{\mathop {\mathrm {Re}}\nolimits}
\renewcommand{\Im}{\mathop {\mathrm {Im}}\nolimits}
\newcommand{\sgn}{\mathop {\mathrm {sgn}}\nolimits}

\def\bfa{\mathbf a}
\def\bfb{\mathbf b}
\def\bfc{\mathbf c}
\def\bfd{\mathbf d}
\def\bfe{\mathbf e}
\def\bff{\mathbf f}
\def\bfg{\mathbf g}
\def\bfh{\mathbf h}
\def\bfi{\mathbf i}
\def\bfj{\mathbf j}
\def\bfk{\mathbf k}
\def\bfl{\mathbf l}
\def\bfm{\mathbf m}
\def\bfn{\mathbf n}
\def\bfo{\mathbf o}
\def\bfp{\mathbf p}
\def\bfq{\mathbf q}
\def\bfr{\mathbf r}
\def\bfs{\mathbf s}
\def\bft{\mathbf t}
\def\bfu{\mathbf u}
\def\bfv{\mathbf v}
\def\bfw{\mathbf w}
\def\bfx{\mathbf x}
\def\bfy{\mathbf y}
\def\bfz{\mathbf z}

\def\bfA{\mathbf A}
\def\bfB{\mathbf B}
\def\bfC{\mathbf C}
\def\bfD{\mathbf D}
\def\bfE{\mathbf E}
\def\bfF{\mathbf F}
\def\bfG{\mathbf G}
\def\bfH{\mathbf H}
\def\bfI{\mathbf I}
\def\bfJ{\mathbf J}
\def\bfK{\mathbf K}
\def\bfL{\mathbf L}
\def\bfM{\mathbf M}
\def\bfN{\mathbf N}
\def\bfO{\mathbf O}
\def\bfP{\mathbf P}
\def\bfQ{\mathbf Q}
\def\bfR{\mathbf R}
\def\bfS{\mathbf S}
\def\bfT{\mathbf T}
\def\bfU{\mathbf U}
\def\bfV{\mathbf V}
\def\bfW{\mathbf W}
\def\bfX{\mathbf X}
\def\bfY{\mathbf Y}
\def\bfZ{\mathbf Z}

\def\frD{\mathfrak D}
\def\frL{\mathfrak L}

\def\bfw{\mathbf w}
%%% END MATHBF
%%%%%%%%%%%%%%%%%%%%%%%%%%%%%%%
%%%%%%%%%%%%%%%%%%%%%%%%%%%%%%%%%
%%% BEGIN MATHBB

\def\R {{\mathbb R }}
 \def\C {{\mathbb C }}
  \def\Z{{\mathbb Z}}
  \def\H{{\mathbb H}}
\def\K{{\mathbb K}}
\def\N{{\mathbb N}}
\def\Q{{\mathbb Q}}
\def\A{{\mathbb A}}

\def\T{\mathbb T}
\def\P{\mathbb P}

\def\G{\mathbb G}

\def\cD{\EuScript D}
\def\cL{\mathscr L}
\def\cK{\EuScript K}
\def\cM{\EuScript M}
\def\cN{\EuScript N}
\def\cR{\EuScript R}
\def\cW{\EuScript W}
\def\cY{\EuScript Y}
\def\cF{\EuScript F}

\def\bbA{\mathbb A}
\def\bbB{\mathbb B}
\def\bbD{\mathbb D}
\def\bbE{\mathbb E}
\def\bbF{\mathbb F}
\def\bbG{\mathbb G}
\def\bbI{\mathbb I}
\def\bbJ{\mathbb J}
\def\bbL{\mathbb L}
\def\bbM{\mathbb M}
\def\bbN{\mathbb N}
\def\bbO{\mathbb O}
\def\bbP{\mathbb P}
\def\bbQ{\mathbb Q}
\def\bbS{\mathbb S}
\def\bbT{\mathbb T}
\def\bbU{\mathbb U}
\def\bbV{\mathbb V}
\def\bbW{\mathbb W}
\def\bbX{\mathbb X}
\def\bbY{\mathbb Y}

\def\kappa{\varkappa}
\def\epsilon{\varepsilon}
\def\phi{\varphi}
\def\le{\leqslant}
\def\ge{\geqslant}

\def\B{\mathrm B}

\def\la{\langle}
\def\ra{\rangle}

\def\F{{}_2F_1}
\def\FF{{}_2^{\vphantom \C}F_1^\C}

\newcommand{\dd}[1]{\,d\,{\overline{\overline{#1}}} }

\def\lambdA{{\boldsymbol{\lambda}}}
\def\alphA{{\boldsymbol{\alpha}}}
\def\betA{{\boldsymbol{\beta}}}
\def\mU{{\boldsymbol{\mu}}}

\def\1{\boldsymbol{1}}
\def\2{\boldsymbol{2}}

\def\SL{\mathrm{SL}}
\def\even{\mathrm{even}}

\begin{center}
	\Large\bf
	An analog of the Dougall formula
	\\
	 and of the de Branges--Wilson integral
	 
	\bigskip
	 
	\large\sc Yury A. Neretin%
	\footnote{Supported by the grants FWF, P28421, P31591.}
\end{center}

{\small
We derive a beta-integral over $\Z\times \R$, which is a counterpart of the
Dougall  $_5H_5$-formula and of the de Branges--Wilson integral, our integral includes
$_{10}H_{10}$-summation. For a derivation we use a two-dimensional integral transform
related to  representations of the Lorentz group,
this transform is a counterpart of the  Olevskii index transform (a synonym: Jacobi transform).}

\section{The statement}

\COUNTERS

{\bf \punct Gamma function of the complex field.}
Denote by $\Lambda_\C$ the set of all pairs $a|a'\in\C^2$
such that $a-a'\in \Z$. For nonzero $z\in\C$ we denote
$$
z^{a|a'}:=z^a\,\ov z^{\,a'}:=|z|^a \Bigl(\frac{\ov z} z\Bigr)^{a'-a}.
$$
Denote by $\Lambda\subset\Lambda_\C$ the set of all $a|a'\in\Lambda_\C$
satisfying the additional condition: $a+a'$ is pure imaginary. We have
$$
\Bigl|z^{a|a'}\Bigr|=1\qquad \text {for $a|a'\in\Lambda$.}
$$
Elements of $\Lambda$ can be represented as
$$
a|a'=a|-\ov a=\tfrac12(k+is)\Bigr| \tfrac12(-k+is), \qquad \text{where $k\in\Z$, $s\in \R$}.
$$
Let $z$ be a complex variable. Denote the Lebesgue measure by
$$
\dd z:=d\Re z\, d\Im z.
$$
Following \cite{GGR}, define the gamma function of the complex field by
\begin{multline}
\Gamma^\C(a|a'):= 
\frac 1\pi
\int_\C z^{a-1|a'-1} e^{2i\Re z} \dd{z}= \\=
i^{a-a'}
\frac{\Gamma(a)}{\Gamma(1-a')}=
i^{a'-a} \frac{\Gamma(a')}{\Gamma(1-a)}=
\frac{i^{a'-a}}\pi \Gamma(a)\Gamma(a')\sin \pi a'
.
\label{eq:gamma}
\end{multline}

Here $a|a'\in \Lambda_\C\simeq \Z\times \C$.
The  $\Gamma^\C$-function  has poles  at points
$a|a'=-k|-l$, where $k$, $l\in\N$.

\sm

{\bf \punct The statement.}
We derive the following beta-integral:

\begin{theorem}
	\label{th:}
	Let $a_1$, $a_2$, $a_3$, $a_4$ satisfy the conditions
	$a_\alpha>0$, $\sum a_\alpha<1$. Then
\begin{multline}
\frac 1{4\pi^2}
\sum_{k=-\infty}^\infty
\int_{-\infty}^{\infty}
\biggl| (k+is)
\prod_{\alpha=1}^4 \Gamma^\C\Bigl(a_\alpha+\tfrac{k+is}2\bigl|
a_\alpha+\tfrac{-k+is}2
\Bigr)\biggr|^2 \,ds
=\\=
\frac{\prod_{1\le\alpha<\beta\le 4}
	\Gamma^\C(a_\alpha+a_\beta|a_\alpha+a_\beta)}
{\Gamma^\C\bigl(a_1+a_2+a_3+a_4\bigl|a_1+a_2+a_3+a_4\bigr)}.
\label{eq:main}
\end{multline}	
We also can write the left hand side as
\begin{multline}
\frac 1{4\pi^2}
\sum_{k=-\infty}^\infty
\int_{-\infty}^{\infty}
 (k+is)(k-is)\times\\\times
\prod_{\alpha=1}^4 \Gamma^\C\Bigl(a_\alpha+\tfrac{k+is}2\bigl|
a_\alpha+\tfrac{-k+is}2
\Bigr) \,
\Gamma^\C\Bigl(a_\alpha+\tfrac{k-is}2\bigl|
a_\alpha+\tfrac{-k-is}2
\Bigr)
\,ds.
\label{eq:bis}
\end{multline}
Then the identity with the same right-hand side holds for 
\begin{equation}\Re a_\alpha>0,\qquad
\Re \sum a_\alpha<1.
\label{eq:aaa}
\end{equation}
\end{theorem}

{\bf\punct The de Branges--Wilson integral and the Dougall formula.}
Recall that the de Branges--Wilson integral is given by
\begin{equation}
\frac 1{4\pi}
\int_{-\infty}^\infty \left|\frac{\prod_{\alpha=1}^{4}
\Gamma(a_\alpha+is)}{\Gamma(2is)}\right|^2\,ds
=\frac{\prod_{1\le\alpha<\beta\le 4} 
	\Gamma(a_\alpha+a_\beta)}{\Gamma(a_1+a_2+a_3+a_4)}.
\label{eq:wilson}
\end{equation}
This formula was obtained by de Branges  \cite{dB0},\cite{dB} in 1972, a proof was not published;
 the formula was rediscovered by Wilson
\cite{Wil}, 1980; see also \cite{AAR}.
The Dougall $_5H_5$ formula is
\begin{equation}
\sum_{k=-\infty}^\infty
\frac{k+\theta}
{\prod_{\alpha=1}^{4}\Gamma(b_\alpha+\theta+k) \Gamma(b_\alpha-\theta-k)}
=\frac{\sin 2\pi\theta}{2\pi}
\,
\frac{\Gamma(b_1+b_2+b_3+b_4-3)}{\prod_{1\le\alpha<\beta\le 4} \Gamma(b_\alpha+b_\beta-1)}
\label{eq:dougall}
\end{equation}
Setting $b_4+\theta=1$, we get a series $\sum_{k\ge0}$
of type $_5F_4[\dots;1]$, this result was contained in a family of 
identities obtained by
Dougall \cite{Dou}, 1906. It seems that the general bilateral
formula was obtained by Bailey \cite{Bai}, 1936;
see also \cite{AAR}.

Let us explain a similarity of such formulas.
Denote by $I(s)$ the integrand in
(\ref{eq:wilson}) and extend it into a complex domain
writing
$$
|\Gamma(c+is)|^2=\Gamma(c+is)\Gamma(c-is).
$$
Next, consider the sum 
$\sum_{k=-\infty}^\infty I(i\theta+ik)$.
We have
$$\Gamma(2is)\Gamma(-2is)\Bigr|_{s=i(\theta+k)}=
\frac{-4\pi(\theta+k)}{\sin 2\pi\theta},
$$
and we get the left-hand side of (\ref{eq:dougall}) with $b_\alpha=1-a_\alpha$.
Also, we get the same right-hand sides.

In \cite{Ner-jacobi} (see formula (1.28))  there was derived a one-dimensional hybrid of 
(\ref{eq:wilson}) and (\ref{eq:dougall}) including 
 both integration over the real axis and a summation over a lattice
 on the imaginary axis.

Our integral (\ref{eq:main}) can be obtained by formal replacing
$\Gamma$-factors in the integrand (\ref{eq:wilson}) by similar $\Gamma^\C$-factors. The function $\Gamma_\C$
satisfies the reflection identity
\begin{equation}
\Gamma^\C(a|a')\,\Gamma^\C(1-a|1-a')=(-1)^{a'-a}.
\label{eq:reflection}
\end{equation}
Therefore
$$
\Gamma^\C(k+is|-k+is)\, \Gamma^\C(-k-is|k-is)=\frac 1{(k+is)(k-is)}
$$
and we come to the integrand in (\ref{eq:main}).

\sm 

{\bf\punct Further structure of the paper.}
In Section 2 we derive
our integral (\ref{eq:main}).
For a calculation we use a unitary integral transform
$J_{a,b}$
defined in \cite{MN}, see below Subsect. \ref{ss:index}. This transform
 is an analog of a classical
integral transform known under the names 
{\it generalized Mehler--Fock transform, Olevskii transform, Jacobi transform,}
see \cite{Koo-paley}, \cite{Koo}, \cite{Ner-w}.

We write an appropriate family of functions $H_\mu$, and our integral
(\ref{eq:main}) arises
as the identity $\la J_{a,b} H_\mu, J_{a,b} H_\nu\ra=\la  H_\mu,  H_\nu\ra$.

Section 3 contains a further discussion of Theorem \ref{th:}.

\section{Calculation}

\COUNTERS

{\bf \punct Convergence of the integral.}

\begin{lemma}
	Let $a|a'\in\Lambda_\C$, $\lambda|\lambda'\in \Lambda$, i.e.,
	$\lambda'=-\ov\lambda$. Then
	$$
	\Gamma^\C(a+\lambda|a'-\ov\lambda)
	\sim i^{a-a'+\lambda+\ov\lambda} \lambda^{-\tfrac12+a|-\tfrac12+a'}\cdot
	\frac{\lambda^\lambda}{{\ov \lambda}^{\ov\lambda}}
	\cdot
	e^{-\lambda+\ov\lambda},\qquad
	\text{as $|\lambda|\to\infty$.}
	$$
	In particular,
	$$
\bigl|	\Gamma^\C(a+\lambda|a'-\ov\lambda)\bigr|
\sim|\lambda|^{-1+\Re(a+a')}
,\qquad
\text{as $|\lambda|\to\infty$.}
	$$
\end{lemma}

{\sc Remark.} Our expression is single valued. Indeed,
$$
\lambda^\lambda{\ov \lambda}^{\,-\ov\lambda}=
\exp\bigl\{\lambda (\ln\lambda+2\pi i N)- \ov\lambda (\ln\ov\lambda-2\pi i N) \bigr\}
=\lambda^\lambda{\ov \lambda}^{\,-\ov\lambda} \exp\bigl\{2\pi iN(\lambda+\ov\lambda)\bigr\}.
$$
But $\lambda+\ov \lambda\in\Z$, and therefore the result does not depend on a choice
of a branch of $\ln\lambda$.
\hfill $\boxtimes$

\sm

{\sc Proof.} We use two expressions for $\Gamma^\C(a+\lambda|a'-\ov\lambda)$,
namely,
$$
i^{a-a'+\lambda+\ov\lambda}\frac{\Gamma(a+\lambda)}{\Gamma(1-a'+\ov\lambda)}=
i^{-a+a'-\lambda-\ov\lambda}\frac{\Gamma(a'-\ov\lambda)}{\Gamma(1-a-\lambda)}.
$$
If $|\arg \lambda|<\pi-\epsilon$, then we apply the Stirling formula
(see e.g., \cite{AAR}) to the first expression. If $|\arg (-\lambda)|<\pi-\epsilon$,
we apply it to the second expression.
\hfill $\square$

\begin{corollary}
	\label{cor:}
	If the parameters $a_\alpha$ satisfy {\rm(\ref{eq:aaa})}, then the integral
	{\rm(\ref{eq:bis})} absolutely converges.
\end{corollary}

{\bf \punct The Gauss hypergeometric function of the complex field.}
For $h|h'\in\Lambda_\C$ we denote
$$
[h|h']=\tfrac12 \Re(h+h')
.$$
Following \cite{GGR} (see also \cite{GGV}, Subsect. II.3.7) we define the beta function
$B^\C[\cdot]$ and the Gauss hypergeometric
function $\FF[\cdot]$ of the complex field. Let $a|a'$,
$b|b'\in\Lambda$. Then
\begin{equation}
B^\C(a|a',b|b'):=\frac 1\pi \int_\C t^{a-1|a'-1}(1-t)^{b-1|b'-1} \,\dd{t}
=\frac{\Gamma^\C(a|a')\,\Gamma^\C(b|b')}
{\Gamma^\C(a+b|a'+b')}
.
\label{eq:beta}
\end{equation}
The integral absolutely converges iff
$$
[a|a']>0,\quad [b|b']>0,\quad [a|a']+[b|b']<2.
$$
 The right hand side
gives a meromorphic continuation of $B^\C$
 to the whole $\Lambda_\C^2\simeq\Z^2\times\C^2$.

For $a|a'$,
$b|b'$, $c|c'\in\Lambda$ we define the hypergeometric function
\begin{multline}
\FF\Bigl[\begin{matrix}a|a',b|b'\\c|c'\end{matrix};z\Bigr]
:=\\:=
\frac1{\pi B^\C(b|b',c-b|c'-b')}
\int_\C
t^{b-1|b'-1}(1-t)^{c-b-1|c'-b'-1}(1-zt)^{-a|-a'}
\dd{t}.
\label{eq:def-FF}
\end{multline}
The integral has an open domain of convergence 
on any connected component of the set of parameters
$\Lambda_\C^3\simeq \Z^3\times\C^3$, it admits a meromorphic
continuation to the whole set $\Lambda_\C^3$, see \cite{MN},  Section 3.

The functions $\FF[\cdot]$ admit explicit expressions
in terms of sums of products of Gauss hypergeometric functions
$\F$. The standard properties of Gauss hypergeometric functions
can be transformed to similar properties of functions
$\FF$, see \cite{MN}, Section 3.

Below we need the following analog of the Gauss formula for $\F[a,b;c;1]$:
\begin{equation}
\FF\Bigl[\begin{matrix}a|a',b|b'\\c|c'\end{matrix};z\Bigr]
=\frac{\Gamma^\C(c|c')\, \Gamma^\C(c-a-b|c'-a'-b')}
{\Gamma^\C(c-a|c'-a')\,\Gamma^\C(c-b|c'-b')},
\label{eq:Gauss}
\end{equation}
which is valid if
\begin{equation}
[c|c']>[a|a']+[b|b'],
\label{eq:cond}
\end{equation}
see \cite{MN}, Proposition 3.2, the last condition coincides with a condition
of continuity of $\FF[\dots;z]$ at $z=1$.

\sm

{\bf \punct The index hypergeometric transform.%
\label{ss:index}}
Fix real $a$, $b$ such that
$$
0\le a\le1,\quad 0\le b\le 1,\qquad
(a,b)\ne (\pm1, \pm1),\, (\mp1,\pm 1).
$$
Consider the measure on $\C$ given by
$$
\rho_{a,b}(z)\,\dd z=
|z|^{2a+2b-2}|1-z|^{2a-2b}\,\dd{z}.
$$
Let
$$
\lambda|\lambda'=\tfrac{k+is}{2}\Bigl| \tfrac{-k+is}{2}\in\Lambda.
$$
Consider the following function on $\Lambda\simeq \Z\times \R$:
$$
\kappa_{a,b} (\lambda|\lambda')=\kappa_{a,b}(k,s):
= \Bigl|
\lambda\,
\Gamma^\C(a-\lambda|a+\ov\lambda')\,
\Gamma^\C(b+\lambda|b-\lambda')\Bigr|^2
$$
and consider the space $L^2_\even(\Lambda, \kappa_{a,b})$
of even functions on $\Lambda$ with inner product
$$
\la\Phi,\Psi\ra_{L^2_\even(\Lambda, \kappa_{a,b})}:=
\sum_{k=0}^\infty \int_{-\infty}^\infty
\Phi(k,s)\,\ov {\Psi(k,s)}
\kappa_{a,b}(k,s)\,ds.
$$
Next, define the kernel 
on $\C\times \Lambda$  by
\begin{equation}
\cK_{a,b}(\lambda|\lambda')=
\frac1{\Gamma^\C(a+b|a+b)}
\,
\FF\left[\begin{matrix}
a+\lambda|a-\lambda',\,
a-\lambda|a+\lambda'
\\ a+b|a+b
\end{matrix}
;z
\right].
\label{eq:cK}
\end{equation}
In \cite{MN} there was obtained the following statement:

\sm

{\it The operator $J_{a,b}$ defined by 
$$
J_{a,b}f(\lambda|\lambda'):=\int_{\C} \cK(z;\lambda|\lambda')\,f(z)\,\rho_{a,b}(z)\dd z
$$
is a unitary operator}
$$
L^2(\C,\rho_{a,b})\to L^2_\even(\Lambda, \kappa_{a,b}).
$$

\sm

{\bf\punct Application of the Mellin transform.}
We define a Mellin transform $\cM$ on $\C$ as the Fourier
transform on the multiplicative group $\C^\times$ of $\C$.
Since $\C^\times\simeq (\R/2\pi\Z)\times \R$, the Mellin transform
 is reduced to the usual Fourier
transform and Fourier series.  We have
\begin{equation}
\cM f(\xi|\xi')=
\cM f\bigl(\tfrac{l+\tau}2\bigl| \tfrac{-l+\tau}2\bigr)
:=\int_\C t^{\xi-1|\xi'-1}
f(t)\,\dd t,
\label{eq:mellin}
\end{equation}
where $\xi|\xi'\in\Lambda_\C$.
In the cases discussed below a function $f$ on $\C\setminus 0$
is differentiable except the  point $t=1$, where
a singularity has a form $C_1+C_2(1-t)^{h|h'}$, $[h|h']>-1$.
Also in our cases 
the integral (\ref{eq:mellin}) absolutely converges
for $\sigma$ being in a certain strip
$A<[\xi|\xi']<B$, therefore
the Mellin transform is holomorphic
 in the strip.
The inversion formula is given by
$$
f(t)=\frac 1{4\pi^2 i}\sum_{l=-\infty}^{\infty}
\int_{\gamma-i\infty}^{\gamma+i\infty}
t^{-(l+\tau)/2|-(-l+\tau)/2}
\,
 \cM f\bigl(\tfrac{l+\tau}2\bigl| \tfrac{-l+\tau}2\bigr)\, d\tau,
$$
the integration is taken over arbitrary line $\Re \sigma=\gamma$, where $A<\gamma<B$.
We understand the  integral (which can be non absolutely convergent) as
$$
\sum_{l=-\infty}^{\infty}
\int_{\gamma-i\infty}^{\gamma+i\infty}:=
\lim_{N\to\infty,\, R\to\infty}\,\,\, \sum_{-N}^N \int_{\gamma-iR}^{\gamma+iR},\qquad\qquad
$$
The inversion formula holds at points of differentiability of
$f$, also it holds at points of singularities
of the form $C_1+C_2(1-t)^{h|h'}$ with $\Re(h+h')>-0$.
In this case we can  repeat the standard proof of pointwise
inversion formula for the one-dimensional Fourier transform and pointwise convergence
of Fourier series, see, e.~g., \cite{KF}, Sect. VIII.1, VIII.3;
for advanced multi-dimensional versions of the Dini condition,
 see, e.g., \cite{Zhi}, Sect. 9.

%Denote by $L^1(\C^\times)$ the space $L^1$ on $\C$ with respect
%to the measure $|z|^{-2}\dd z$.
%If $f\in L^1(\C^\times)$, then its Mellin transform 
%is continuous.
% For $f_1$, $f_2\in L^1(\C^\times)$
%have a well defined  
Convolution $f_1*f_2$ on $\C^\times$ is defined by
$$
f_1*f_2(t):=\int_{\C} f_1(t/z)f_2(z)\,|z|^{-2}\,\dd z.
$$ 
As usual, we have
$$
\cM(f_1*f_2)=\cM(f_1) \cM(f_2),
$$
this identity holds in intersection of strips of holomorphy
$\cM(f_1)$ and  $\cM(f_2)$.
We also define a function $f^\star(t)=f(t^{-1})$.
Then 
$$
\cM f^\star (\xi|\xi')=\cM f (-\xi|-\xi').
$$
So we have the following corollary of the convolution formula:
\begin{multline}
\int_\C f_1(t)\, f_2(t) \,|t|^{-2}\,\dd t=
\cM (f_1^\star*f_2) (1)=
\\=
\frac 1{4\pi^2 i}\sum_{l=-\infty}^{\infty}
\int_{\gamma-i\infty}^{\gamma+i\infty}
\cM f_1\bigl(-\tfrac{l+\tau}2\bigl|- \tfrac{-l+\tau}2\bigr)
\,
\cM f_2\bigl(\tfrac{l+\tau}2\bigl| \tfrac{-l+\tau}2\bigr)
\, d\sigma,
\label{eq:ff}
\end{multline}
where the integration contour is contained in the intersection
of domains of holomorphy of $\cM f_1^\star$ and $\cM f_2$.

\begin{lemma}
	\label{l:1}
	{\rm a)}
	Let $[q|q']<0$.
Then the Mellin transform sends a function 
$$t^{p|p'}(1-t)^{q|q'}$$ to
$$
\Gamma^\C(q+1|q'+1)
\frac{\Gamma^\C\bigl(p+\xi\bigl|p'+\xi')}
{\Gamma^\C\bigl(p+q+1+\xi\bigl|p'+q'+1+\xi')},
$$
it is holomorphic in the strip
$$
-[p|p']<[\xi|\xi']<-[p|p']-[q|q'].
$$

{\rm b)}
Assume that
$$
[c|c']-[a|a']-[b|b']>-1.
$$
Then the
Mellin transform sends a function 
$$
\FF\left[\begin{matrix}
a|a',b|b'\\c|c'
\end{matrix};t  \right]
$$
to the function
\begin{equation}
\frac{
\Gamma^\C(a-\xi|a'-\xi')\,
\Gamma^\C(b-\xi|b'-\xi')\,
\Gamma^\C(\xi|\xi')\,
\Gamma^\C(1-c+\xi|1-c'+\xi')
}{\Gamma^\C(a|a')\,\Gamma^\C(b|b')\,\Gamma^\C(1-c|1-c')}
\label{eq:mellin-F}
\end{equation}
defined in the strip
%\begin{equation}
%\max\bigl(-2, \Re(-2+c+c')\bigr)<\Re(\xi+\xi')<\min(\Re(a+a'), \Re(b+b')).
%\label{eq:strip-F}
%\end{equation}
\begin{equation}
\max\bigl(-1, -1+[c|c']\bigr)<[\xi|\xi']<\min\bigl([a|a'], [b|b]\bigr)
\label{eq:strip-F}
\end{equation}
{\rm(}if it is non-empty{\rm)}
\end{lemma}

{\sc Proof.} The statement a) follows from the definition of $B^\C$-function (\ref{eq:beta}).

\sm

b) The function $\FF[z]$ has singularities at
$z=0$, $1$, $\infty$
with asymptotics of the form
\begin{align}
\FF[z]&\sim A_1+A_2 |z|^{1-c|1-c'}\qquad&\text{as $z\to 0$;}
\label{eq:as1}
\\
\FF[z]&\sim B_1+B_2 (1-z)^{c-a-b|c'-a'-b'} \qquad&\text{as $z\to 1$;}
\label{eq:as2}
\\
\FF[z]&\sim C_1 z^{-a|-a'}+ C_2 z^{-b|-b'} \qquad&\text{as $z\to \infty$,}
\label{eq:as3}
\end{align}
see \cite{MN}, Theorem 3.9. This gives us a strip of convergence
of the Mellin transform.

We must evaluate the integral
\begin{equation}
%\frac 1{\pi B^\C(b|b',c-b|c'-b')}
\int_{\C} z^{\xi-1|\xi'-1} \int_\C t^{b-1|b'-1} 
(1-t)^{c-b-1|c'-b'-1}(1-zt)^{-a|-a'}\,\dd t\,\dd z
.
\label{eq:double}
\end{equation}
We change an order of the integration and  integrate in $z$:
$$
\int_{\C} z^{\xi-1|\xi'-1} (1-zt)^{-a|-a'}\,\dd z= 
\pi B^\C(\xi|\xi', -a+1|-a'+1) t^{-\xi|-\xi'}.
$$ 
Integrating in $t$ we again met a $B^\C$-function
and after simple cancellations and an application
of the reflection formula (\ref{eq:reflection})
we come to (\ref{eq:mellin-F}). 
The successive integration is valid under conditions
(\ref{eq:strip-F}).

We must justify the change of order of integrations.
In fact, (\ref{eq:double}) is absolutely convergent as a double
integral, i.e., 
$$
\int_{\C} \int_\C  |z|^{\Re(\xi+\xi')-2} t^{b+b'-2} 
|1-t|^{\Re(c+c'-b-b')-2}|1-zt|^{-\Re(a+a')}\,\dd t\,\dd z
<\infty.
$$
It is a special case of integral
(\ref{eq:double}), we integrate it 	successively
 in $z$ and in $t$ under the same condition as for successive
 integration in (\ref{eq:double}).
 \hfill $\square$

\sm

\begin{lemma}
	\label{l:2}
	Let $\lambda|\lambda'\in \Lambda$ and
	\begin{equation}
	a>0, \quad b>0,\quad \mu>0,\quad a+\mu<1,\quad b+\mu<1.
	\label{eq:convergence}
	\end{equation}
	Then
	\begin{multline}	
	\int_{\C} z^{a+b-1|a+b-1} (1-z)^{-b-\mu|-b-\mu}
	\FF\left[\begin{matrix}
	a+\lambda|a-\lambda', a-\lambda|a+\lambda'
	\\
	a+b|a+b
	\end{matrix};z\right]\,\dd z
	=\\=
	\frac{\Gamma^\C(a+b|a+b) \,
\Gamma^\C(\mu+\lambda|\mu+\lambda')\,\Gamma^\C(\mu-\lambda|\mu-\lambda')	}
{\Gamma^\C(a+\mu|a+\mu) \, \Gamma^\C(b+\mu|b+\mu)}
\label{eq:aux}
	\end{multline}
\end{lemma}
 
 {\sc Proof.} We apply formula (\ref{eq:ff}) assuming
 $$f_1:=z^{a+b-1|a+b-1} (1-z)^{-b-\mu|-b-\mu},\qquad
 f_2:=\FF[\dots;z].
 $$
We evaluate Mellin transforms of $f_1$, $f_2$  applying  Lemma \ref{l:1}. 
In the integrand in the right
 hand side of (\ref{eq:ff}) we get a product of two factors. The first
 factor is 
 \begin{equation*}
 \frac{\Gamma^\C(-b-\mu+1|-b-\mu+1)\,\Gamma^\C(a+b-\xi|a+b-\xi')}
 {\Gamma^\C(a-\mu+1+\xi|a-\mu+1+\xi')},
 \end{equation*}
 it is holomorphic un the strip
 \begin{equation}
 -a-b<[\xi|\xi']<\mu-a.
 \label{eq:strip1}
 \end{equation}
 The second factor
 \begin{equation*}
 \frac{\Gamma^\C(a+b|a+b)\,\Gamma^\C(a+\lambda-\xi|a-\lambda-\xi')
 	\Gamma^\C(a-\lambda-\xi|a+\lambda-\xi')
 	\Gamma^\C(\xi|\xi')(-1)^{\xi-\xi'}}
 {\Gamma^\C(a+\lambda|a+\lambda') \Gamma^\C(a-\lambda|a-\lambda')\Gamma^\C(a+b-\xi|a+b-\xi')}
 \end{equation*}
 is holomorphic in the strip
 \begin{equation}
 a+b-1<[\xi|\xi']<a.
 \label{eq:strip2}
 \end{equation}
 It $a$, $b$ are sufficiently small, then
 strips  (\ref{eq:strip1}) and (\ref{eq:strip2})
 have a non-empty intersection and we can apply formula
 (\ref{eq:ff}).
 Two factors $\Gamma^\C(a+b-\xi|a+b-\xi')$ cancel and we get 
 a factor independent on $\xi$ and the integral
 \begin{multline*}
 \sum%_{l=-\infty}
 \int
 \frac{\Gamma^\C(a+\lambda-\xi|a-\lambda-\xi')
 \Gamma^\C(a-\lambda-\xi|a+\lambda-\xi')
 \Gamma^\C(\xi|\xi')(-1)^{\xi-\xi'}}{\Gamma^\C(a-\mu+1+\xi|a-\mu+1+\xi')}\,
d\tau.
 \end{multline*}
 The integrand  up to a constant factor is a Mellin transform
 of a function $\FF[\dots;z]$. By the inversion formula,  integral
 (\ref{eq:aux}) converts to
 \begin{equation*}
 \frac{\Gamma^\C(a+b|a+b)\, 
 	\Gamma^\C(-b-\mu-1|-b-\mu-1)}{\Gamma^\C(a-\mu-1|a-\mu-1)}
 \,\FF\left[\begin{matrix}
 a+\lambda|a-\lambda',a-\lambda|a+\lambda'\\
 a-\mu+1|a-\mu+1
 \end{matrix};1\right].
 \end{equation*} 
For sufficiently small $a$ we can apply the Gauss identity (\ref{eq:Gauss}).
Thus we get (\ref{eq:aux}) for sufficiently small $a$, $b>0$.
	
	Keeping in mind (\ref{eq:as1})--(\ref{eq:as3}),
	we can easily verify that the  integral 
	 in the left-hand side of (\ref{eq:aux})
	converges for
	$$
\Re	a>0, \quad \Re b>0,\quad \Re \mu>0,\quad \Re(a+\mu)<1,\quad \Re(b+\mu)<1.
	$$
Thus, under these conditions the left hand side is holomorphic.	
The right hand side also is holomorphic in this domain. Therefore,
they coincide.
 \hfill $\square$
 
 \sm

{\bf\punct Proof of Theorem 1.1.}
Let $\in \R$.
Consider a function $H_\mu(z)$ on $\C$ given by
$$
H_\mu(z):=(1-z)^{-a-\mu|-a-\mu}.
$$

\begin{lemma}
	{\rm a)} $H_\mu\in L^2(\C, \rho_{a,b})$ iff
	$$
0<	2\mu<1-a-b.
	$$
	
	{\rm b)}
	$$
	\la H_\mu,H_\nu\ra_{L^2(\C,\rho_{a,b})}=
	\frac{\Gamma^\C(a+b|a+b)\,\Gamma^\C(\mu+\nu|\mu+\nu)}
	{\Gamma^\C(a+b+\mu+\nu|a+b+\mu+\nu)}
	.$$	
\end{lemma}

The statement a) is trivial, b) is reduced to $B^\C$-function.
\hfill $\square$

\sm

The $J_{a,b}$-image of $H_\mu$ is done by Lemma \ref{l:2}.
Since $J_{a,b}$ is unitary, we have
$$
\la H_\mu,H_\nu\ra_{L^2(\C,\rho_{a,b})}
=\la J_{a,b}H_\mu,J_{a,b}H_\nu\ra_{L^2_\even(\Lambda,\kappa_{a,b})}
.
$$
This is Theorem \ref{th:}, where the parameters $a_1$, $a_2$, $a_3$, $a_4$ 
are $a$, $b$, $\mu$, $\nu$. Our calculation is valid for positive 
reals
$a_1$, $a_2$, $a_3$, $a_4$ satisfying conditions 
$a_1+a+2+2a_3<1$, $a_1+a+2+2a_4<1$.
For extending the identity to the domain (\ref{eq:aaa})
we refer to Corollary \ref{cor:}, the integral (\ref{eq:bis})
is holomorphic in the domain (\ref{eq:aaa}), the right hand side also is holomorphic.

\section{Final remarks}

\COUNTERS

{\bf \punct  Barnes--Ismagilov integrals.}
Let $p\le q$.
Let $a_\alpha|a_\alpha'$, $b_\alpha|b_\alpha'\in\Lambda$. Following Ismagilov \cite{Ism},
we  define integrals of the form
\begin{multline*}
I_{p,q}[a,b;z]:=
\frac 1{2\pi i}\sum_{k=-\infty}^\infty
\int\limits_{-i\infty}^{i\infty}
\prod_{\alpha=1}^{p}
\Gamma^\C\Bigl(a_\alpha+\tfrac{k+\sigma}2\Bigl|a'_\alpha+\tfrac{-k+\sigma}2
\Bigr)
\times\\\times
\prod_{\beta=1}^q 
\Gamma^\C\Bigl(b_\beta+\tfrac{-k-\sigma}2\Bigl|b_\beta'+\tfrac{k-\sigma}2
\Bigr) z^{(k+\sigma)/2\bigl|(-k+\sigma)/2}\,d\sigma.
\end{multline*}
By \cite{Ism}, Lemma 2,
such integrals admit a representation of the form
$$
\sum_{j=1}^q \gamma_j(\cdot) _p^{\vphantom \C}F_{q-1}[\cdot, z]\, _pF_{q-1}[\cdot, z],
$$
where $\gamma_j(\cdot)$ are products of $\Gamma$-factors and parameters
of hypergeometric functions $_p^{\vphantom \C}F_{q-1}[\cdot, z]$ are linear expressions in $a_\alpha$, $b_\beta$
and $a'_\alpha$,
$b'_\beta$. 
It is reasonable to claim that integrals
$I_{p,q}=:\,_p^{\vphantom \C}F_{q-1}^\C$ are hypergeometric functions of the complex field.

By Lemma \ref{l:2}.b, the functions $\FF$ defined by (\ref{eq:def-FF}) are compatible
with this definition. 
Ismagilov considered $_4^{\vphantom \C}F^\C_3[\dots;1]$-expressions that are counterparts of the Racah
coefficients for unitary representations of the Lorentz group  $\SL(2,\C)$. 
Our theorem is an example of a hypergeometric identity
for $_5^{\vphantom \C} F^\C_4[\dots;1]$.

\sm

Some integrals of  type (\ref{eq:main}) with products of $\Gamma^\C$-functions were obtained by Kells \cite{Kel} and Derkachov, Manashov, and Valinevich in \cite{DMV}. 

\sm

{\bf \punct A difference problem.} The de Branges--Wilson integral, the Dougall formula,
and our integral (\ref{eq:main}) are representatives of beta integrals
in the sense of Askey \cite{Ask}.
  Quite often integrands $w(x)$ in beta integrals are 
weight functions for systems of hypergeometric orthogonal polynomials.
In particular, orthogonal polynomials corresponding to the de Branges--Wilson integral
are the Wilson polynomials, see \cite{Wil}, \cite{AAR}. Recall that they are even eigenfunctions of the following difference operator
 \begin{equation*}
L f(s)=
\frac{\prod\limits_{\alpha=1}^4 (a_\alpha+is)}{2is(1+2is)}
\bigl(f(s-i)-f(s))
+
\frac{\prod\limits_{\alpha=1}^4 (a_\alpha-is)}{-2is(1-2is)}
\bigl(f(s+i)-f(s)),
\end{equation*}
where $i^2=-1$.
If an integrand $w(x)$ of a beta integral  decreases as a power function, then
only finite number of moments $\int x^n w(x)\,dx$ converge;
however in this case a beta integral can be a weight
for a finite system of hypergeometric orthogonal polynomials
(this phenomenon was firstly observed by Romanovski in
\cite{Rom}), the system of orthogonal polynomials 
related to the Dougall $_5H_5$ formula
was obtained in \cite{Ner-w}. On the other hand, such finite systems
 are discrete parts of spectra of explicitly solvable Sturm--Liouville
 problems (see, e.g., \cite{Gro}, \cite{Ner-two}). 
 
 In the case of our integral (\ref{eq:main}) the integrand decreases
 as $|k+is|^{2\sum a_j -8}$, we have no orthogonal polynomials.
 However a difference Sturm--Liouville problem can be formulated. We consider a space of
 meromorphic even functions $\Phi(\lambda|\lambda')$ on $\Lambda_\C$,
 a weight on $\Lambda\subset \Lambda_\C$ defined by the integrand
 (\ref{eq:main}), and the following commuting difference operators:
 \begin{multline*}
 \frL \Phi(\lambda|\lambda')=
 \frac{\prod_{\alpha=1}^4 (a_\alpha+\lambda)}{2\lambda(1+2\lambda)}
 \bigl(\Phi(\lambda+1|\lambda')-\Phi(\lambda|\lambda'))
 +\\+
 \frac{\prod_{\alpha=1}^4 (a_\alpha-\lambda)}{-2\lambda(1-2\lambda)}
 \bigl(\Phi(\lambda-1|\lambda')-\Phi(\lambda|\lambda'));
 \end{multline*}  
  \begin{multline*}
 \frL \Phi(\lambda|\lambda')=
 \frac{\prod_{\alpha=1}^4 (a_\alpha+\lambda')}{2\lambda'(1+2\lambda')}
 \bigl(\Phi(\lambda|\lambda'-1)-\Phi(\lambda|\lambda'))
 +\\+
 \frac{\prod_{\alpha=1}^4 (a_\alpha-\lambda')}{-2\lambda'(1-2\lambda')}
 \bigl(\Phi(\lambda|\lambda'+1)-\Phi(\lambda|\lambda')).
 \end{multline*} 
See  a simpler pair of difference operators of this kind
in \cite{MN}. On the other hand, see a one-dimensional operator with
continuous spectrum similar to
$L$ 
in Groenevelt \cite{Gro}.

 \tt
\noindent
Yury Neretin\\
Math. Dept., University of Vienna/c.o Pauli Institute \\
\&Institute for Theoretical and Experimental Physics (Moscow); \\
\&MechMath Dept., Moscow State University;\\
\&Institute for Information Transmission Problems;\\
URL: http://mat.univie.ac.at/$\sim$neretin/	

\end{document}